\documentclass{arxiv}

\def\cprime{$'$}
\newcommand{\HF}{Hartree-Fock model}
\newcommand{\HFT}{Hartree-Fock model with temperature}
\newcommand{\FE}{free energy}
\usepackage{mathrsfs}
\usepackage{dsfont}
\newcommand\R{{\ensuremath {\mathbb R} }}

\newcommand\N{{\ensuremath {\mathbb N} }}

\renewcommand\phi{\varphi}

\newcommand{\gH}{\mathfrak{H}}
\newcommand{\gS}{\mathfrak{S}}

\newcommand{\wto}{\rightharpoonup}

\newcommand{\cF}{\mathcal F}

\newcommand{\cM}{\mathcal M}

\newcommand{\tr}{{\rm tr}\,}
\newcommand{\cE}{\mathcal{E}}

\newcommand{\cK}{\mathcal{K}}

\newcommand\ii{{\ensuremath {\infty}}}
\newcommand\pscal[1]{{\ensuremath{\left\langle #1 \right\rangle}}}
\newcommand{\norm}[1]{\left\|\, #1 \,\right\|}

\date{February 7, 2008}
\begin{document}

%%%%%%%%%%%%%%%%%%%%%%%%%%%%%%%%%%%%%%%%%%%%%%%%%%%%%%%%%%%%%%%%%%%%%%%%%%%%
%%%%%%%%%%%%%%%%%%%%%%%%%%%%%%%%%%%%%%%%%%%%%%%%%%%%%%%%%%%%%%%%%%%%%%%%%%%%

\markboth{J. Dolbeault, P. Felmer \& M. Lewin}{Hartree-Fock model with temperature}

%%%%%%%%%%%%%%%%%%% Publisher's Area please ignore %%%%%%%%%%%%%%%%%%%%%%%%%
%
%\catchline{}{}{}{}{}
%
%%%%%%%%%%%%%%%%%%%%%%%%%%%%%%%%%%%%%%%%%%%%%%%%%%%%%%%%%%%%%%%%%%%%%%%%%%%%

\title{Stability of the Hartree-Fock model with temperature}

\author{Jean Dolbeault}

\address{Ceremade (UMR CNRS no. 7534), Universit\'e Paris-Dauphine, Place de Lattre de Tassigny, F-75775 Paris C\'edex 16, France. E-mail: \email{\textsf{dolbeaul@ceremade.dauphine.fr}}}

\author{Patricio Felmer}

\address{Universidad de Chile, Facultad de Ciencias Fisicas y Matem\'aticas, Depto Ingenier\'{\i}a Matem\'atica, Blanco Encalada 2120, Piso 5, Santiago, Chile. E-mail: \email{\textsf{pfelmer@dim.uchile.cl}}}

\author{Mathieu Lewin}

\address{D\'epartement de Math\'ematiques (UMR CNRS no. 8088), Universit\'e de Cergy-Pontoise, Site de Saint-Martin, 2, avenue Adolphe Chauvin, 95302 Cergy-Pontoise Cedex France.\newline E-mail: \email{\textsf{Mathieu.Lewin@math.cnrs.fr}}}

\maketitle

\begin{history}
%\received{(Day Month Year)}
%\revised{(Day Month Year)}
%\accepted{(Day Month Year)}
%\comby{(xxxxxxxxxx)}
\end{history}

\begin{abstract} {\sc Abstract.} This paper is devoted to the \HFT~in the euclidean space. For large classes of free energy functionals, minimizers are obtained as long as the total charge of the system does not exceed a threshold which depends on the temperature. The usual \HF~is recovered in the zero temperature limit. An orbital stability result for the Cauchy problem is deduced from the variational approach. \end{abstract}



% 26-XX REAL FUNCTIONS [See also 54C30]
% 26Dxx Inequalities {For maximal function inequalities, see 42B25; for functional inequalities, see 39B72; for probabilistic inequalities, see 60E15 } [See also 46E27, 60Bxx]
% 26D15 Inequalities for sums, series and integrals
%
% 35-XX PARTIAL DIFFERENTIAL EQUATIONS
% 35Jxx Partial differential equations of elliptic type [See also 58J10, 58J20]
% 35J10 Schroedinger operator [See also 35Pxx]
% 35Qxx Equations of mathematical physics and other areas of application [See also 35J05, 35J10, 35K05, 35L05]
% 35Q40 Equations from quantum mechanics
%
% 47-XX OPERATOR THEORY
% 47Bxx Special classes of linear operators
% 47B34 Kernel operators
% 47Gxx Integral, integro-differential, and pseudo-differential operators [See also 58Jxx]
% 47G20 Integro-differential operators [See also 34K30, 35R10, 45J05, 45K05]
%
% 81-XX QUANTUM THEORY
% 81Qxx General mathematical topics and methods in quantum theory
% 81Q10 Selfadjoint operator theory in quantum theory, including spectral analysis
% 81Vxx Applications to specific physical systems
% 81V45 Atomic physics
%
% 82-XX STATISTICAL MECHANICS, STRUCTURE OF MATTER
% 82Bxx Equilibrium statistical mechanics
% 82B10 Quantum equilibrium statistical mechanics (general)

%%%%%%%%%%%%%%%%%%%%%%%%%%%%%%%%%%%%%%%%%%%%%%%%%%%%%%%%%%%%%%%%%%%%%%%%%%%%
%%%%%%%%%%%%%%%%%%%%%%%%%%%%%%%%%%%%%%%%%%%%%%%%%%%%%%%%%%%%%%%%%%%%%%%%%%%%
\section{Introduction}

The \emph{\HFT\/} is a simple extension of the \HF~\cite{Lieb-Simon-77}. The minimization of the \emph{\FE\/} determines an equilibrium state in the presence of a thermal noise, whose effect is to populate electronic states corresponding to excited energy levels. Compared to the energy functional of the standard \HF, the \FE~of the \HFT~$T$ has an additional term, which is the product of $T$ with an entropy term. The standard \HF~corresponds to the limit case $T=0\,$.

The \HFT~has been studied \cite{Bach-Lieb-Solovej-94,Lions-88} in the case of an entropy based on the function $\beta(\nu)=\nu\log\nu+(1-\nu)\log(1-\nu)\,$. The main drawback in this case is that the \FE~is unbounded from below when the model is considered in the euclidean space.

It has been recently established \cite{Dolbeault-Felmer-Loss-Paturel-05} that the \FE~for a quantum system described by a \emph{mixed state\/} or a density operator, in the presence of an external potential $V$, can be uniformly bounded from below by a functional depending on $V$. This follows from interpolation estimates of Gagliardo-Nirenberg type for systems, which are equivalent to Lieb-Thirring inequalities. The result of course depends on the convex function on which the entropy is based. We shall not make a direct use of such inequalities, but our approach will rely on classes of entropies for which the \FE~is well defined and semi-bounded if $V(x)=Z/|x|\,$.

We shall add another ingredient, which comes from the modeling in quantum mechanics. We will impose that \emph{occupation numbers\/} take values in $[0,1]$ because the electrons are fermions. This can be formalized by considering  entropies based on convex functions with support in $[0,1]$ and by assuming that they take infinite values on $\mathbb R\setminus[0,1]\,$. In this paper, we shall adopt the language of \emph{density operators\/}: a state of the system is represented by a trace-class self-adjoint operator $\gamma$ satisfying $0\leq\gamma\leq1\,$. Occupation numbers $n_i$ are the eigenvalues of the density operator: \hbox{$\gamma=\sum_{i\geq1}n_i\,|u_i\rangle\langle u_i|\,$}. When only a finite number of eigenvalues is non-zero, the operator is of finite rank. We shall say that it describes a \emph{pure state\/} when all non-zero occupation numbers are equal to~$1$, \emph{i.e.\/} $\gamma$ is an orthogonal projector. Minimizers of the \HF~are of this form, but this will not be the case in general when the temperature $T>0$ is big enough: in this case the corresponding minimizer is a \emph{mixed state,\/} that has some non-integer eigenvalues.

A practical consequence of the positivity of the occupation numbers is that the exchange term is dominated by the Poisson term in the Hartree-Fock energy functional. In the case of a bounded domain, not much is changed when the temperature is turned on. A pure Hartree model (without the exchange term), which is usually called the Schr\"odinger-Poisson system, has already been studied \cite{MR1889607,DFM}.

The purpose of this paper is to deal with the whole euclidean space and we shall see that many families of free energies can be considered, each of them giving rise to minimizers under some conditions which eventually depend on the temperature. This is not surprising. Similar phenomena have indeed been observed in other models of mechanics, for instance in kinetic theory. The most studied example is probably the case of self-gravitating gases in astrophysics. In such models, relaxation mechanisms are not so well known, but stationary states can be observed and the corresponding \FE~functional, which is also called the energy-Casimir functional, can then be recovered quite easily. Stationary solutions being characterized as minimizers of the \FE, orbital stability appears as a simple consequence for the collisionless version of the equation \cite{MR1991989}. In presence of collisions, the picture is compatible with hydrodynamic or diffusive limits, which provides a very natural link with nonlinear diffusion equations \cite{MR2338354}. Similar results are known in other frameworks like fluid mechanics \cite{MR0180051} or nonlinear Schr\"odinger equations \cite{MR677997}. The striking observation is that many different stationary states can be used simultaneously for measuring the orbital stability of a given solution to the Cauchy problem.

In this paper, we will focus on the mathematical aspects of this model. We consider \emph{free energy} functionals of the form
$$\gamma\mapsto \cE^{\rm HF}(\gamma)-TS(\gamma)$$
where $\cE^{\rm HF}$ is the Hartree-Fock energy and the entropy takes the form 
$$S(\gamma):=-\,\tr\big(\beta(\gamma)\big)\,,$$
$\beta$ being a convex function on $[0,1]\,$. We study under which conditions (for instance on the total number of electrons, or global charge) such a free energy has a minimizer.

Next we consider the stability of the associated time-dependent equation describing the evolution of the state $\gamma\,$,
\begin{equation}
i\,\frac{d\gamma}{dt}=[H_\gamma,\gamma]\;,
\label{time-dependent_intro}
\end{equation}
which is called the \emph{von Neumann equation.\/} Here $H_\gamma$ is a self-adjoint operator depending on $\gamma\,$ but not on the function $\beta$. Orbital stability is then a straightforward consequence of the variational method and the existence theory for the Cauchy problem \cite{MR0413843,MR0411439}. Since we consider a large class of entropies based on various functions $\beta\,$, we thereby construct a large class of associated stable states of the time-dependent equation \eqref{time-dependent_intro}.

The paper is organized as follows. In the next section, we state our main results and only give the shortest proofs. Other proofs are detailed in Section~\ref{sec:proofs}.

%%%%%%%%%%%%%%%%%%%%%%%%%%%%%%%%%%%%%%%%%%%%%%%%%%%%%%%%%%%%%%%%%%%%%%%%%%%%
%%%%%%%%%%%%%%%%%%%%%%%%%%%%%%%%%%%%%%%%%%%%%%%%%%%%%%%%%%%%%%%%%%%%%%%%%%%%
\section{Results}

%%%%%%%%%%%%%%%%%%%%%%%%%%%%%%%%%%%%%%%%%%%%%%%%%%%%%%%%%%%%%%%%%%%%%%%%%%%%
\subsection{Definition of the free energy}

The Hartree-Fock energy, written in terms of the density matrix $\gamma\,$, reads
\begin{equation}
\cE^{\rm HF}(\gamma):=\tr\big((-\Delta)\,\gamma\big) -Z\int_{\R^3}\frac{\rho_\gamma}{|x|}\;dx + \frac12\, D(\rho_\gamma,\rho_\gamma) - \frac{1}{2}\iint_{\R^3\times\R^3}\frac{|\gamma(x,y)|^2}{|x-y|}\;dx\,dy\;.\label{def_HF_energy}
\end{equation}
The first term of the right hand side of \eqref{def_HF_energy} is the kinetic energy of the electrons and the second is the electrostatic interaction with a classical nucleus of charge $Z$, located at $0\in\R^3$. We use the notation
$$D(f,g):=\iint_{\R^3\times\R^3}\frac{f(x)\,g(y)}{|x-y|}\;dx\,dy$$
for the classical Coulomb interaction between two densities of charge $f$ and $g\,$. Such a quantity is well-defined as soon as, for instance, $f,g\in L^{6/5}(\R^3)\,$, by the Hardy-Littlewood-Sobolev inequality \cite{Lieb-Loss-01}. The last two terms of the r.h.s. of \eqref{def_HF_energy} describe the interaction between the electrons. The classical electrostatic interaction $D(\rho_\gamma,\rho_\gamma)/2$ is usually called the \emph{direct term}. The last term of \eqref{def_HF_energy} is a purely quantum term, which is a consequence of the Pauli exclusion principle, called the \emph{exchange term}. The Hartree-Fock energy \eqref{def_HF_energy} was first studied from a mathematical point of view by Lieb and Simon \cite{Lieb-Simon-77} and then by various other authors \cite{Lions-87,Bach-Lieb-Solovej-94,Bach-Lieb-Loss-Solovej-94}. 

The \HFT~\cite{Bach-Lieb-Solovej-94} relies on the free energy, which is defined as
\begin{equation}
\cE_Z^\beta(\gamma):=\cE^{\rm HF}(\gamma)+T\,\tr\big(\beta(\gamma)\big)\;.
\label{def_free_energy}
\end{equation}
The function $\beta$ in \eqref{def_HF_energy} is a convex function defined on $[0,1]$ with values in $\R\,$. To take into account the constraint on the occupation numbers, we extend it to $+\infty$ on $\R\setminus[0,1]\,$. The quantity $-\,\tr\big(\beta(\gamma)\big)$ is the entropy and $T$ is the temperature. In quantum mechanics, a common choice is $\beta (\nu)=\nu\log \nu+(1-\nu)\log(1-\nu)\,$, but the model makes sense only on a compact subset of $\R^3$ or for finite-rank density matrices \cite{Lions-88}. We shall investigate other entropies \cite{Dolbeault-Felmer-Loss-Paturel-05} for which the problem can be studied in the whole euclidean space.

\medskip In the whole paper, we assume that
\begin{description}
\item[(A1)] $\beta$ is a strictly convex $C^1$ function on $(0,1)\,$,\medskip
\item[(A2)] $\beta(0)=0$ and $\beta\geq0$ on $[0,1]\,$.
\end{description}
Then we introduce a modified Legendre transform of $\beta$
\begin{equation*}
g(\lambda):=\mathop{\rm argmin}_{0\leq\nu\leq1}\left(\lambda\,\nu+\beta(\nu)\right)\,,
\end{equation*}
\emph{i.e.}
$$g(\lambda)=\sup\Big\{\inf\big\{(\beta')^{-1}(-\lambda),\,1\big\}\,,\;0\Big\}\;.$$
Notice that $g$ is a nonincreasing function with $0\leq g\leq 1\,$. Also define
$$\beta^*(\lambda):=\lambda\,g(\lambda)+(\beta\circ g)(\lambda)\;.$$

As we will be interested in the minimization of the energy \eqref{def_HF_energy} under a charge constraint of the form $\tr\gamma=q\,$, we can substract $T\,\beta'(0)\,\tr\gamma$ to the definition of the free energy \eqref{def_free_energy}. Hence we shall assume from now on that
\begin{description}
\item[(A3)] $\beta$ is a nonnegative $C^1$ function on $[0,1)$ and $\displaystyle \beta'(0)=0\,$,
\end{description}
without much loss of generality. With this assumption, $g_{|\R^+}\equiv 0$ and $g$ is positive on $(-\infty,0)\,$. Then we also assume that $\beta^*$ satisfies the following condition,
\begin{description}
\item[(A4)] $\displaystyle \sum_{j\geq1}\,j^2\,\big|\,\beta^*\!\left(-Z^2/(4\,T\,j^2)\right)\big|<\ii\,$.
\end{description}
The role of this $T$-dependent condition is to ensure that the ground state free energy is finite, using an estimate based on the eigenvalues of the no-spin hydrogen atom. We recall that the eigenvalues of $-\Delta-Z/|x|$ are $-Z^2/(4\,j^2)$ and that each eigenvalue is degenerated with multiplicity~$j^2\,$.

%---------------------------------------------------------------------------
\begin{example}\label{Example:Power}A typical example \cite{Dolbeault-Felmer-Loss-Paturel-05} is
\begin{equation*}
\beta(\nu)=\nu^m
\label{example_power}
\end{equation*}
which satisfies \textnormal{\textbf{(A1)--(A4)}} as long as $1<m<3\,$. In this special case
$$g(\lambda)=\left\{
\begin{array}{ll}
\min\left\{\left({\textstyle \frac{-\lambda}{m}}\right)^{\frac1{m-1}},1\right\} \quad& \text{if}\quad \lambda<0\,,\vspace*{12pt}\\
0 & \text{otherwise}\,,
\end{array}\right.$$
and
$$\beta^*(\lambda)=-(m-1)\left({\textstyle \frac{-\lambda}{m}}\right)^{\frac m{m-1}}\quad\mbox{if}\quad -m<\lambda<0\;.$$
\end{example}
%---------------------------------------------------------------------------

\medskip We define next an adapted functional setting for the definition of the energy in terms of density matrices. Let us introduce the following space of operators
$$\gH:=\left\{ \gamma:L^2(\R^3)\to L^2(\R^3)\;|\;\gamma=\gamma^*,\;\gamma\in\gS_1\,,\;\sqrt{-\Delta}\,|\gamma|\,\sqrt{-\Delta}\in\gS_1\right\}\,,$$
where we have denoted by $\gS_1$ the space of trace-class operators \cite{Reed-Simon,Simon}. This is a Banach space when endowed with the norm
$$\norm{\gamma}_\gH=\tr\,|\gamma|\,+\tr\big(\sqrt{-\Delta}\,|\gamma|\,\sqrt{-\Delta}\,\big)\;.$$
We also introduce the convex closed subset of $\gH$ defined by
$$\cK:=\left\{ \gamma\in\gH\;|\;0\leq\gamma\leq1\right\}\;.$$

Since $\beta$ is convex and $\beta(0)=0\,$, we have $0\leq\beta(\nu)\leq \beta(1)\,\nu$ on $[0,1]\,$. Hence for any $\gamma\geq0\,$, we have $0\leq \beta(\gamma)\leq \beta(1)\,\gamma\,$. This proves that $\beta(\gamma)\in\gS_1$ when $\gamma\in\gS_1$ and therefore $\tr\big(\beta(\gamma)\big)$ is well-defined for any $\gamma\in\cK\,$.

For any nonnegative operator $\gamma\,$, we use the shorthand notation
$$\tr\big((-\Delta)\,\gamma\big):=\tr\left(\sqrt{-\Delta}\,\gamma\,\sqrt{-\Delta}\,\right)\in\R\cup\{+\ii\}\;.$$
Of course $\tr\!\left((-\Delta)\,\gamma\right)$ is finite when $\gamma\in\cK\,$.

For any $\gamma\in \cK\,$, we define the associated density of charge:
$$\rho_\gamma(x)=\gamma(x,x)\in L^1(\R^3)$$
where $\gamma(x,y)$ is the kernel of the trace-class operator $\gamma\,$. Using the spectral decomposition of $\gamma\,$, the following classical inequality is easily proved:
\begin{equation*}
\forall\,\gamma\in\cK\;,\qquad \norm{\nabla\sqrt{\rho_\gamma}}_{L^2(\R^3)}^2\leq \tr\big(\sqrt{-\Delta}\,\gamma\,\sqrt{-\Delta}\,\big)\;.
\label{estim_rho_nabla}
\end{equation*}
Hence we have
\begin{equation*}
\forall\,\gamma\in\cK\;,\qquad \norm{\sqrt{\rho_\gamma}}_{H^1(\R^3)}^2\leq \norm{\gamma}_\gH\;.
\label{estim_rho_H1}
\end{equation*}
and, as a consequence, $\rho_\gamma\in L^1(\R^3)\cap L^3(\R^3)\subset L^{6/5}(\R^3)\,$. This also shows that $\int_{\R^3}\frac{\rho_\gamma(x)}{|x|}\,dx<\ii$ and $D(\rho_\gamma,\rho_\gamma)<\ii\,$, by the Hardy-Littlewood-Sobolev inequality~\cite{Lieb-Loss-01}. Finally, from \hbox{$\gamma\geq0\,$}, it follows by the Cauchy-Schwarz inequality for sequences that
\begin{equation}
|\gamma(x,y)|^2\leq \rho_\gamma(x)\,\rho_\gamma(y)\quad \text{for a.e. } (x,y)\in\R^3\times\R^3,
\label{estim_exchange_term}
\end{equation}
so that
$$\iint_{\R^3\times\R^3}\frac{|\gamma(x,y)|^2}{|x-y|}\;dx\,dy\leq D(\rho_\gamma,\rho_\gamma)<\ii\;.$$

These considerations show that $\cE_Z^\beta$ introduced in \eqref{def_HF_energy} is a well-defined funtional on $\cK\,$. We are also interested in minimizing $\cE_Z^\beta$ under a constraint corresponding to a closed convex subset of $\cK\,$,
$$\cK_q:=\{\gamma\in\cK\;|\;\tr\gamma=q\}\;.$$

%%%%%%%%%%%%%%%%%%%%%%%%%%%%%%%%%%%%%%%%%%%%%%%%%%%%%%%%%%%%%%%%%%%%%%%%%%%%
\subsection{The linear case}\label{sec_linear}

We now recall the properties of the linear case corresponding to
\begin{equation}
 \cF_Z^\beta(\gamma):=\tr\big((-\Delta)\,\gamma\big) -Z\int_{\R^3}\frac{\rho_\gamma}{|x|}\;dx +T\,\tr\big(\beta(\gamma)\big)\;.
\label{def_linear_energy}
\end{equation}
This corresponds to the case where the last two terms in \eqref{def_HF_energy} are removed, \emph{i.e.\/}
$$\cE_Z^\beta(\gamma)=\cF_Z^\beta(\gamma)+\frac12\, D(\rho_\gamma,\rho_\gamma) - \frac{1}{2}\iint_{\R^3\times\R^3}\frac{|\gamma(x,y)|^2}{|x-y|}\;dx\,dy\;.$$
Assume that $\beta$ satisfies \textnormal{\textbf{(A1)--(A4)}}. Then a straightforward minimization \cite{Dolbeault-Felmer-Loss-Paturel-05} gives\begin{equation}
\inf_{ \gamma\in\cK}\cF_Z^\beta(\gamma) = \sum_{j\geq1}\,j^2\,\beta^*\Big(\frac{\lambda_j}T\Big) > -\ii\;, \label{value_linear}
\end{equation}
where $\lambda_j=-\frac{Z^2}{4\,j^2}$ are the negative eigenvalues of $-\Delta-Z/|x|$ with multiplicity $j^2\,$. A ground state for \eqref{value_linear} is formally given by 
$$\gamma=g\left({\textstyle\frac 1T\,\big(-\Delta-\frac Z{|x|}\big)}\right)\;.$$
However, this state does not necessarily belong to $\cK$. Although its kinetic energy is finite by \eqref{value_linear}, its trace
\begin{equation}
q_{\rm max}^{\rm lin}(T):=\tr\!\left(g\left({\textstyle\frac 1T\,\big(-\Delta-\frac Z{|x|}\big)}\right)\right)=\sum_{j\geq1}\,j^2\,g\Big(\frac{\lambda_j}T\Big)\in(0,\ii]
\label{def_q_max}
\end{equation}
could in principle be infinite. Then the minimization problem with constraint
\begin{equation}
\inf_{ \gamma\in\cK_q}\cF_Z^\beta(\gamma)\label{value_linear_constraint}
\end{equation}
admits a minimizer for all $q\geq0$ if $q_{\rm max}^{\rm lin}=\ii\,$, whereas it has a minimizer if and only if $q\in[0,q_{\rm max}^{\rm lin}]$ if $q_{\rm max}^{\rm lin}<\ii\,$. In all cases, this minimizer solves the equation
$$\gamma=g\left({\textstyle\frac 1T(-\Delta-\frac Z{|x|}-\mu)}\right)$$
for some $\mu\leq0\,$, a Lagrange multiplier which is chosen to ensure that the condition $\tr\gamma=q$ is satisfied. The function $\mu\mapsto q(\mu):=\tr(g\big((-\Delta-Z/|x|-\mu)/T\big))$ is indeed non decreasing and satisfies $q(\mu)=0$ for $\mu<0$, $|\mu|$ large enough and $q(\mu)\to q_{\rm max}^{\rm lin}$ when $\mu\to0\,$. The range of the function $q(\mu)$ gives all possible $q$'s for which there is a minimizer for \eqref{value_linear_constraint}.

When $q_{\rm max}^{\rm lin}$ is finite, it can take very small values, depending on the temperature~$T\,$. For instance, let us take $\beta(\nu)=\nu^m$ as in Example~\ref{Example:Power}. In this case we see that $q_{\rm max}^{\rm lin}<\ii$ if and only if $m<5/3\,$, as summarized in Table \ref{fig:summary_linear}. For $T>Z^2/(4\,m)\,$,
\begin{equation}
q_{\rm max}^{\rm lin}(T)=\left(\frac{Z^2}{4\,T\,m}\right)^{\frac1{m-1}}\sum_{j\geq1}\, j^{2-\frac 2{m-1}}=\left(\frac{Z^2}{4\,T\,m}\right)^{\frac1{m-1}}\!\zeta\left(2\,\frac{m-2}{m-1}\right)
\label{formula_q_max_example}
\end{equation}
where $\zeta$ denotes the Riemann zeta function. We observe that $q_{\rm max}^{\rm lin}(T)\to0$ as $T\to\ii\,$. We shall observe a similar result in the nonlinear case.
%---------------------------------------------------------------------------
\begin{table}\small
\centering
\begin{tabular}{|c|c|c|}
\hline
$1<m<5/3$ & $5/3\leq m<3$ & $m\geq3$\\
\hline\hline&&\\
$q_{\rm max}^{\rm lin}(T)<\ii$ & $q_{\rm max}^{\rm lin}(T)=\ii$ &Linear\\
& & energy \eqref{def_linear_energy}\\
Existence iff& Existence & unbounded\\
$0\leq q\leq q_{\rm max}^{\rm lin}(T)$ & $\forall\, q\geq0$ & from\\
& & below\\
\cline{1-2}
\multicolumn{2}{|c|}{Energy \eqref{def_linear_energy} bounded from below}& \\
\hline
\end{tabular}
\caption{Existence and non-existence of minimizers with a finite trace for \eqref{value_linear_constraint} for $\beta(\nu)=\nu^m$.}
\label{fig:summary_linear}
\end{table}
%---------------------------------------------------------------------------

%%%%%%%%%%%%%%%%%%%%%%%%%%%%%%%%%%%%%%%%%%%%%%%%%%%%%%%%%%%%%%%%%%%%%%%%%%%%
\subsection{Minimization of the free energy}\label{Sec:Minimization}

As a consequence of \eqref{estim_exchange_term} and \eqref{value_linear}, we obtain that $\cE_Z^\beta$ is bounded from below on~$\cK$ for any $\beta$ satisfying \textbf{(A1)--(A4)}. Hence we can define
\begin{equation}
I^\beta_Z(q):=\inf\left\{\cE_Z^\beta(\gamma)\;|\;\gamma\in\cK\;\text{and}\;\tr(\gamma)=q\right\}
\label{def_min_charge_constraint}
\end{equation}
and
\begin{equation*}
I^\beta_Z:=\inf\left\{\cE_Z^\beta(\gamma)\;|\;\gamma\in\cK\right\}=\inf_{q\geq0} I^\beta(q)\;.
\label{def_global_min}
\end{equation*}
It is easily seen that $\cE_Z^\beta$ is continuous for the topology of $\gH\,$. If we do not put any external potential, that is if we take $Z=0\,$, $I^\beta_0(q)$ can be computed explicitly:
%---------------------------------------------------------------------------
\begin{lemma}[Ground state energy with $Z=0$]\label{lem:value_free} Assume that $\beta$ satisfies \textnormal{\textbf{(A1)--(A4)}} for some $T>0$. Then we have for any $q\geq0$
$$I^\beta_0(q)=\inf_{\substack{0\leq n_i\leq 1\\ \sum_{i\geq1} n_i=q}}\beta(n_i)=\beta'(0)\,q=0\;.$$
\end{lemma}
%---------------------------------------------------------------------------
\begin{proof} Let $Z=0\,$. By \eqref{estim_exchange_term} and using $\tr\big((-\Delta)\,\gamma\big)\geq0\,$, we have that
$$I^\beta_0(q) \geq \inf_{\substack{\gamma\in\cK\\ \tr\gamma=q}} \tr\big(\beta(\gamma)\big)= \inf_{\substack{0\leq n_i\leq 1\\ \sum_{i\geq1} n_i=q}}\sum_{i\geq1}\beta(n_i)=0\;.$$
Now let $\epsilon>0$ and $\gamma=\sum_{1\leq i\leq K}\nu_i\,|\phi_i\rangle\langle\phi_i|$ be a finite-rank operator such that
$$\tr\big(\beta(\gamma)\big)\leq \inf_{\substack{\gamma'\in\cK\\ \tr\gamma'=q}} \tr\big(\beta(\gamma')\big)+\epsilon=\epsilon\;.$$
As the $\phi_i$'s can be chosen arbitrarily, we can assume that $\phi_i\in H^2(\R^3)$ for all $1\leq i\leq K\,$, in which case $(-\Delta)\,\gamma\in\gS_1\,$. Let $U_\eta:L^2(\R^3)\to L^2(\R^3)$ be the dilatation unitary operator defined as $U_\eta(\phi)(x):=\eta^{3/2}\phi(\eta\,x)$ and which is such that $U_\eta^*=U_{1/\eta}\,$. Notice that $U_\eta^*\,\gamma\,U_\eta\in\cK$ and that $\tr\big(\beta(U_\eta^*\,\gamma\,U_\eta\big)=\tr\big(\beta(\gamma)\big)$ for all $\eta>0\,$. Using the equality $(-\Delta)\,U_\eta=\eta^2\,U_\eta(-\Delta)\,$, we infer
$$\tr\big((-\Delta)\,U_\eta^*\,\gamma\,U_\eta\big)=\frac1{\eta^2}\,\tr\big((-\Delta)\,\gamma\big)\;.$$
Similarly, the kernel of $U^*_\eta\,\gamma\,U_\eta$ is $(U^*_\eta\,\gamma\,U_\eta)(x,y)=\eta^{-3}\,\gamma(x/\eta,y/\eta)\,$. Hence $\rho_{\,U^*_\eta\gamma U_\eta}(x)=\eta^{-3}\,\rho_\gamma(x/\eta)$ and
$$D(\rho_{\,U^*_\eta\gamma U_\eta},\rho_{\,U^*_\eta\gamma U_\eta})=\frac 1\eta\,D(\rho_\gamma,\rho_\gamma)\;,$$
$$\iint_{\R^6}\frac{|(U^*_\eta\,\gamma\,U_\eta)(x,y)|^2}{|x-y|}\;dx\,dy=\frac 1\eta\,\iint_{\R^6}\frac{|\gamma(x,y)|^2}{|x-y|}\;dx\,dy\;.$$
Hence we have
$$\cE_0^\beta(U^*_\eta\,\gamma\,U_\eta)\leq C/\eta + \tr\big(\beta(\gamma)\big)\;.$$
Taking first $\eta\to\ii$ and then $\epsilon\to0$ yields the result. \end{proof}

We now state our main result.
%---------------------------------------------------------------------------
\begin{theorem}[Minimization for the HF model with temperature]\label{thm_exists} Assume that $\beta$ satisfies \textnormal{\textbf{(A1)--(A4)}} for some $T>0\,$.
\begin{enumerate}
\item For every $q\geq0\,$, the following statements are equivalent:
\begin{description}
\item[$(i)$] all minimizing sequences $(\gamma_n)_{n\in\N}$ for $I^\beta_Z(q)$ are precompact in $\cK\,$,
\item[$(ii)$] $I^\beta_Z(q)<I^\beta_Z(q')$ for all $q\,$, $q'$ such that $0\leq q'<q\,$.
\end{description}
\item Any minimizer $\gamma$ of $I^\beta_Z(q)$ satisfies the self-consistent equation
$$\gamma=g\big((H_\gamma-\mu)/T\big)\;,\quad H_\gamma=-\Delta-\frac{Z}{|x|}+\rho_\gamma\ast|\cdot|^{-1}-\frac{\gamma(x,y)}{|x-y|}$$
for some $\mu\leq0\,$.
\item The minimization problem $I^\beta_Z(q)$ has no minimizer if $q\geq 2\,Z+1\,$.
\item Problem $I^\beta_Z$ always has a minimizer $\bar\gamma\,$. It satisfies the self-consistent equation
$$\bar\gamma=g(H_{\bar\gamma}/T)\;.$$
\end{enumerate}
\end{theorem}
%---------------------------------------------------------------------------
\begin{remark} We shall prove below in Lemma~\ref{lem:binding_inequalities} (also see Lemma~\ref{lem:value_free}) that
\begin{equation*}
\forall\,q\,,\;q'\;\mbox{such that}\; 0\leq q'\leq q\;,\quad I^\beta_Z(q)\leq I^\beta_Z(q-q')+I^\beta_0(q')=I^\beta_Z(q-q')\;.
\label{binding_inequalities_rmk}
\end{equation*}
\end{remark}
%---------------------------------------------------------------------------

The proof of Theorem~\ref{thm_exists} is given below in Section~\ref{proof_thm_exists}. Let us investigate the validity of Condition $(ii)$ in Theorem~\ref{thm_exists}. We first give a bound on the largest possible charge $q_{\rm max}^{\rm HF}$ of minimizers by comparing with the linear model.
%---------------------------------------------------------------------------
\begin{proposition}\label{prop_nonexist_g}Assume that $\beta$ satisfies \textnormal{\textbf{(A1)--(A4)}} for some $T>0$. Then the minimization problem \eqref{def_min_charge_constraint} has \emph{no solution} if $q\geq q_{\rm max}^{\rm lin}$ as defined in \eqref{def_q_max}. \end{proposition}
%---------------------------------------------------------------------------
As a consequence of Proposition~\ref{prop_nonexist_g} and Part (3) of Theorem~\ref{thm_exists}, we obtain that the largest possible charge $q_{\rm max}^{\rm HF}$ for the nonlinear problem \eqref{def_min_charge_constraint} satisfies
$$q_{\rm max}^{\rm HF}\leq \min\left\{q_{\rm max}^{\rm lin}\;,\; 2\,Z+1\right\}\;.$$
In case of Example~\ref{Example:Power}, for $\beta(\nu)=\nu^m$, the largest possible charge $q_{\rm max}^{\rm HF}$ converges to zero as $T\to\ii\,$, by \eqref{formula_q_max_example}.

\begin{proof} Let $\gamma$ be a minimizer of $I^\beta_Z(q)$ for some $q\,$. By Theorem~\ref{thm_exists}, it solves the self-consistent equation $\gamma=g\big((H_\gamma-\mu)/T\big)$ for some multiplier $\mu\leq0\,$. Next we notice that
\begin{equation}
\rho_\gamma\ast\frac{1}{|\cdot|}-\frac{\gamma(x,y)}{|x-y|}\geq0
\label{estim_direct_echange}
\end{equation}
in the sense of operators on $L^2(\R^3)\,$. In fact, we notice that
\begin{multline*}
\pscal{\left(\rho_{|\phi\rangle\langle\phi|}\ast\frac{1}{|\cdot|}-\frac{\phi(x)\,\phi(y)}{|x-y|}\right)\psi,\psi} \\
=\iint\frac{\phi(x)^2\,\psi(y)^2-\phi(x)\,\phi(y)\,\psi(x)\,\psi(y)}{|x-y|}\;dx\,dy\geq0
\end{multline*}
by the Cauchy-Schwarz inequality. Then \eqref{estim_direct_echange} follows from the decomposition $\gamma=\sum_{j\geq1}\,n_j\,|\phi_j\rangle\langle\phi_j|$ with $n_j\geq0\,$.
{}From \eqref{estim_direct_echange} we deduce that
$$H_\gamma\geq -\Delta -\frac{Z}{|x|}$$
in the sense of self-adjoint operators on $L^2(\R^3)\,$. Since $g$ is nonincreasing, we infer
$$q=\tr\big(g\big((H_\gamma-\mu)/T\big)\big)\leq \tr\big(g(H_\gamma/T\big)\leq \tr\big(g((-\Delta -Z/|x|)/T)\big)=q_{\rm max}^{\rm lin}\;.$$
\end{proof}

On the other hand we can give a positive lower bound on $q_{\rm max}^{\rm HF}\,$.
%---------------------------------------------------------------------------
\begin{proposition}\label{prop_exists_g} Assume that $\beta$ satisfies \textnormal{\textbf{(A1)--(A4)}} for some $T>0$. Then for all $q$ such that
\begin{equation}
0\leq q\leq \min\Big\{\sum_{j\geq1}g\left({\textstyle\frac{-(Z-q)^2}{4\,T\,j^2}}\right)\, , \, Z\Big\}\,,
\label{estim_q_exists}
\end{equation}
Condition $(ii)$ in Theorem~\ref{thm_exists} is satisfied. \end{proposition}
%---------------------------------------------------------------------------
\begin{remark} With $\lambda_j(\varepsilon)=-\varepsilon ^2/(4\,j^2)$, if $\sum_{j\geq1}g(\lambda_j(\varepsilon)/T)=\ii$ for any $\varepsilon>0\,$, existence of a minimizer for $I^\beta_Z(q)$ holds for all $q\leq Z\,$, whereas, if $\sum_{j\geq1}g(\lambda_j(\varepsilon)/T)$ is finite for some $\varepsilon>0\,$, we get the existence at least on an interval $[0,q_{\rm max}]$ for some $q_{\rm max}>0\,$, because the right hand side of \eqref{estim_q_exists} does not converge to $0$ as $q\to0\,$. It is natural to conjecture that existence holds true if we replace $\sum_{j\geq1}g(\lambda_j(Z-q)/T)$ by $\sum_{j\geq1}j^2\,g(\lambda_j(Z-q)/T$ in the r.h.s. of \eqref{estim_q_exists}, like in the linear case, but we have been unable to prove it.
%---------------------------------------------------------------------------

If $\beta(\nu)=\nu^m$ with $1<m<3\,$, then $\sum_{j\geq1}g(\lambda_j(\varepsilon)/T)<\ii$ for any $\varepsilon>0\,$, and so existence holds true on some interval $[0,q_{\rm max}]\,$. \end{remark}

The proof of Proposition~\ref{prop_exists_g} is given below in Section~\ref{sec:prop_exists_g}. The case $T=0$ is well known \cite{Lieb-Simon-77,Lions-87} and it is not difficult to see that estimates are uniform as $T\to0_+\,$. Summarizing, we have found the following existence result. 
%---------------------------------------------------------------------------
\begin{corollary}[Existence of minimizers for the HF model with temperature]\label{cor_exists} Let $T\geq 0\,$. Assume that $\beta$ satisfies \textnormal{\textbf{(A1)--(A3)}}, and \textnormal{\textbf{(A4)}} if $T$ is positive. Then there exists $q_{\rm max}>0$ such that the minimization problem \eqref{def_min_charge_constraint} has a minimizer for any $q\in[0,q_{\rm max}]\,$.\end{corollary}
%---------------------------------------------------------------------------

%%%%%%%%%%%%%%%%%%%%%%%%%%%%%%%%%%%%%%%%%%%%%%%%%%%%%%%%%%%%%%%%%%%%%%%%%%%%
\subsection{Orbital stability}

An interesting consequence of Theorem~\ref{thm_exists} and Corollary~\ref{cor_exists} is that the set of all minimizers is orbitally stable for the von Neumann time-dependent equation, which reads
\begin{equation}
\left\{\begin{array}{l}
\displaystyle i\,\frac{d\gamma}{dt}=\big[H_\gamma\,,\,\gamma\big]\;,\vspace*{12pt}\\
\gamma(0)=\gamma_0\in\cK\;.
\end{array}\right.
\label{vonNeumann}
\end{equation}
It was proved in \cite{MR0411439} that \eqref{vonNeumann} has a global-in-time solution $t\mapsto\gamma(t)\in C^1(\R,\cK)$ for all fixed $\gamma_0\in\cK_q$ and that $\tr\big(\gamma(t)\big)=q\,$, $\cE^{\rm HF}_Z\big(\gamma(t)\big)=\cE^{\rm HF}_Z(\gamma_0)$ for all $t\in\R\,$. Since $\beta(\gamma)$ commutes with $H_\gamma$, it is clear that $\tr\big(\beta(\gamma(t))\big)$ is also conserved.

For given $\beta$, $Z>0$, $q>0$ and $T\geq0$, let $\cM$ be the set of all minimizers for \eqref{def_min_charge_constraint}. We shall say that $\cM$ is \emph{orbitally stable\/} if and only if for any $\epsilon>0\,$, there exists $\eta>0$ such that for all $\gamma_0\in\cK_q$ with ${\rm dist}(\gamma_0,\cM)\leq\eta\,$, if $t\mapsto\gamma(t)$ is a solution of \eqref{vonNeumann} with initial data $\gamma_0\,$, we have ${\rm dist}(\gamma(t),\cM)\leq\epsilon$ for all $t\in\R\,$. Here ${\rm dist}(\gamma,\cM):=\inf_{\delta\in\cM}\|\gamma-\delta\|_\gH\,$. 

As a consequence of the continuity of $\cF_Z^\beta$ and of the variational approach of Section~\ref{Sec:Minimization}, we have the 
%---------------------------------------------------------------------------
\begin{proposition}[Orbital stability] Assume that $\beta$ satisfies \textnormal{\textbf{(A1)--(A4)}} for some $T>0$ and that $(i)$ holds true in Theorem~\ref{thm_exists} for some $q>0\,$. Then $\cM$ is orbitally stable. \end{proposition}
%---------------------------------------------------------------------------

It is interesting to emphasize that we have been able to construct a large class of orbitaly stable states for the von Neumann equation \eqref{vonNeumann}. Indeed, any function~$\beta$ satisfying the above assumptions gives rise to an orbitaly stable set $\cM$ in $\cK$.

%%%%%%%%%%%%%%%%%%%%%%%%%%%%%%%%%%%%%%%%%%%%%%%%%%%%%%%%%%%%%%%%%%%%%%%%%%%%
%%%%%%%%%%%%%%%%%%%%%%%%%%%%%%%%%%%%%%%%%%%%%%%%%%%%%%%%%%%%%%%%%%%%%%%%%%%%
\section{Proofs of Theorem~\ref{thm_exists} and Proposition~\ref{prop_exists_g}}\label{sec:proofs}

For the sake of simplicity, we shall assume that $T=1$. Proving the same results for any $T>0$ does not add any difficulty.

%%%%%%%%%%%%%%%%%%%%%%%%%%%%%%%%%%%%%%%%%%%%%%%%%%%%%%%%%%%%%%%%%%%%%%%%%%%%
\subsection{Proof of Theorem~\ref{thm_exists}}\label{proof_thm_exists}

We recall that according to [Reed-Simon, Thm VI.26] \cite{Reed-Simon}, $\gS_1$ is the dual of the space of compact operators acting on $L^2(\R^3)\,$. Hence we can endow the Banach space $\gH$ with the weak-$\ast$ topology for which $\gamma_n\wto\gamma\in\gH$ means
$$\tr(\gamma_n\,K)\to\tr(\gamma\,K)\quad \text{and}\quad \tr(\sqrt{-\Delta}\,\gamma_n\sqrt{-\Delta}K)\to\tr(\sqrt{-\Delta}\,\gamma\,\sqrt{-\Delta} K)$$
for all compact operators $K\,$. The convex set $\cK$ is closed for the strong topology of $\gH$ and also closed for this weak-$\ast$ topology. Of course our main problem will be that when $\gamma_n\wto\gamma$ in $\cK\,$, there could be a loss of mass at infinity in such a way that $\tr\gamma<\liminf_{n\to\infty}\tr\gamma_n\,$. Indeed the linear functional $\gamma\mapsto\tr\gamma$ is continuous but not weakly-$\ast$ continuous on $\gH\,$, \emph{i.e.\/} the sets
$$\cK_q:=\{\gamma\in\cK\;|\;\tr\gamma=q\}$$
are not closed for the weak-$\ast$ topology.

The proof of Theorem~\ref{thm_exists} follows some classical ideas which have been introduced in various papers \cite{Lieb-Simon-77,Lions-87,Bach-92,HLS2,HLS3}.

\paragraph*{Step 1. Properties of the energy.}\hfill
%---------------------------------------------------------------------------
\begin{lemma}\label{lem:weakly_continuous_and_coercive} Assume that $\beta$ satisfies \textbf{(A1)--(A4)} with $T=1\,$. Then $\cE_Z^\beta$ is weakly-$\ast$ lower semi-continuous (wlsc-$\ast$) on $\cK\,$. For every $q\geq0\,$, it is coercive on
$$\{\gamma\in\cK\;|\;\tr\gamma\leq q\}=\bigcup_{0\leq q'\leq q}\cK_{q'}\;.$$
\end{lemma}
%---------------------------------------------------------------------------
\begin{proof} Consider a sequence $\{\gamma_n\}_{n\in\N}$ such that $\gamma_n\wto\gamma\in\cK\,$. Fatou's Lemma gives $\tr\big(\sqrt{-\Delta}\,\gamma\,\sqrt{-\Delta}\,\big)\leq\liminf_{n\to\infty}\tr\big(\sqrt{-\Delta}\,\gamma_n\,\sqrt{-\Delta}\,\big).$ As $\{\gamma_n\}_{n\in\N}$ is bounded in~$\cK\,$, $\gamma_n(x,y)$ is bounded in $H^1(\R^3\times\R^3)$ and $\sqrt{\rho_{\gamma_n}}$ is bounded in $H^1(\R^3)\,$. Hence, up to a subsequence, we can assume that $\rho_{\gamma_n}\to\rho_\gamma$ strongly in $L^p_{\rm loc}\,$, $1\leq p<3\,$, and a.e., $\gamma_n(x,y)\to\gamma(x,y)$ strongly in $L^q_{\rm loc}\,$, $2\leq q<12/5\,$, and a.e. Hence \eqref{estim_exchange_term} and Fatou's Lemma give
\begin{multline*}
\iint_{\R^3\times\R^3}\frac{\rho_\gamma(x)\,\rho_\gamma(y)-|\gamma(x,y)|^2}{|x-y|}\;dx\,dy\\
\leq \liminf_{n\to\infty}\left( \iint_{\R^3\times\R^3}\frac{\rho_{\gamma_n}(x)\,\rho_{\gamma_n}(y)-| \gamma_n(x,y)|^2}{|x-y|}\;dx\,dy\right)
\end{multline*}
and
$$\tr\big(\beta(\gamma)\big)\leq\liminf_{n\to\infty}\tr\big(\beta(\gamma_n)\big)\;.$$
Also
$$\lim_{n\to\ii}\int_{\R^3}\frac{\rho_{\gamma_n}(x)}{|x|}\;dx=\int_{\R^3}\frac{\rho_{\gamma_n}(x)}{|x|}\;dx\;.$$
All this shows that $\cE_Z^\beta$ is wlsc-$\ast$ on $\cK\,$. Now we have
\begin{eqnarray*}
\cE_Z^\beta(\gamma) &\geq & \frac12\,\tr\big((-\Delta)\,\gamma\big)+\tr\left(\left(-\frac 12\,\Delta-\frac{Z}{|x|}\right)\,\gamma\right)\\
& \geq & \frac12\,\tr\big((-\Delta)\,\gamma\big) -\,2\,Z^2\,q
\end{eqnarray*}
which proves that $\cE_Z^\beta$ is coercive as stated. \end{proof}

%---------------------------------------------------------------------------
\begin{lemma}\label{lem:conv_equiv_no_mass_at_infinity} Consider a minimizing sequence $\{\gamma_n\}_{n\in\N}\subset\cK$ for $I^\beta_Z(q)$ and assume that $\gamma_n\wto\gamma\in\cK\,$. Then $\gamma_n\to\gamma$ for the strong topology if and only if $\tr\gamma=q\,$. \end{lemma}
%---------------------------------------------------------------------------
\begin{proof} Assume that $\gamma_n\wto\gamma$ and $\tr\gamma=q\,$. Then as $\cE_Z^\beta$ is wlsc-$\ast$ on $\cK\,$, $\gamma$ is a minimizer for $I^\beta_Z(q)\,$. Hence $\lim_{n\to\ii}\cE_Z^\beta(\gamma_n)= \cE_Z^\beta(\gamma)\,$. By the proof of Lemma~\ref{lem:weakly_continuous_and_coercive}, we see that this implies in particular $\lim_{n\to\ii}\tr\big(\sqrt{-\Delta}\,\gamma_n\,\sqrt{-\Delta}\,\big)= \tr\big(\sqrt{-\Delta}\,\gamma\,\sqrt{-\Delta}\,\big)\,$. This is enough \cite{Reed-Simon} to obtain that $\gamma_n\to\gamma$ for the strong topology of $\gH\,$. \end{proof}

\paragraph*{Step 2. Binding Inequalities.}\hfill
%---------------------------------------------------------------------------
\begin{lemma}\label{lem:binding_inequalities} Assume that $\beta$ satisfies \textbf{(A1)--(A4)} with $T=1\,$. For every $q$, $q'$ such that $0\leq q'<q\,$, we have
\begin{equation}
I^\beta_Z(q)\leq I^\beta_Z(q')\;.
\label{binding_large_inequalities}
\end{equation}
If equality holds with $0\leq q'<q\,$, then there exists a minimizing sequence for $I^\beta_Z(q)$ which is not precompact. \end{lemma}
%---------------------------------------------------------------------------
\begin{proof} We consider two states $\gamma\in\cK_{q'}$ and $\gamma'\in\cK_{q-q'}$ such that $\cE_Z^\beta(\gamma)\leq I^\beta_Z(q')+\epsilon$ and $\cE_0^\beta(\gamma')\leq I^\beta_0(q-q')+\epsilon\,$. By density of finite-rank operators in $\cK$ and of $C^\ii_c$ in $L^2(\R^3)\,$, we can assume that $\gamma=\sum_{j=1}^Kn_j\,|\phi_j\rangle\langle\phi_j|$ and $\gamma'=\sum_{j=1}^{K'}n'_j\,|\phi'_j\rangle\langle\phi'_j|$ where the $\phi_j$ and $\phi_j'$ are smooth functions with support in a ball $B(0,R)\subset\R^3\,$. Now we introduce the translation unitary operator $V_\tau$ defined as $(V_\tau\phi)(x)=\phi(x-\tau e)$ where $e$ is a fixed vector in $\R^3\,$. We use the shorthand notation $\gamma'_\tau:=V_\tau^*\,\gamma'\,V_\tau\,$. For $\tau$ large enough, we have $\gamma\,\gamma'_\tau=\gamma'_\tau\,\gamma=0\,$, hence $\gamma+\gamma'_\tau\in\cK_q$ and $\tr\big(\beta(\gamma+\gamma'_\tau)\big)=\tr\big(\beta(\gamma)\big)+\tr\big(\beta(\gamma')\big)\,$. An easy calculation shows that
$$I^\beta_Z(q)\leq \cE_Z^\beta(\gamma+\gamma'_\tau)=\cE_Z^\beta(\gamma)+\cE_0^\beta(\gamma')+O(1/\tau)\;.$$
Taking first $\tau\to\ii$ and then $\epsilon\to0$ yields \eqref{binding_large_inequalities} using $I^\beta_0(q-q')=0$ by Lemma~\ref{lem:value_free}.

If there is an equality in \eqref{binding_large_inequalities} for some $q$ and $0\leq q'<q\,$, a non compact minimizing sequence is built in the same way. \end{proof}

Notice that by Assumption \textbf{(A3)} and Lemma~\ref{lem:value_free}, Lemma~\ref{lem:binding_inequalities} implies that $q\mapsto I^\beta_Z(q)$ is nonpositive, nonincreasing on $[0,\ii)\,$.

\paragraph*{Step 3. Proof that $(i)\Leftrightarrow(ii)\,$.} The implication $(i)\Rightarrow(ii)$ was already proved in Lemma~\ref{lem:binding_inequalities}. We now prove that $(ii)\Rightarrow(i)\,$. To this end, we consider a minimizing sequence $\{\gamma_n\}_{n\in\N}\subset\cK_q\,$. By Lemma~\ref{lem:weakly_continuous_and_coercive}, $\{\gamma_n\}_{n\in\N}$ is bounded in $\gH\,$. We can therefore assume that, up to a subsequence, $\gamma_n\wto\gamma\,$. We argue by contradiction and assume that $\gamma_n$ does not strongly converge to $\gamma\,$. By Lemma~\ref{lem:conv_equiv_no_mass_at_infinity}, this is equivalent to \hbox{$\tr(\gamma)\neq q\,$}. Notice that by Fatou's Lemma, $\tr\gamma\leq\liminf_{n\to\infty}\tr\gamma_n=q\,$, hence there exists a $0\leq q'<q$ such that $\tr\gamma=q'\,$. Modifying $\{\gamma_n\}_{n\in\N}$ if necessary, we can assume that $\{\gamma_n\}_{n\in\N}$ is of finite rank and that $(-\Delta)\,\gamma_n\in\gS_1$ for all $n\,$.

Now we follow a truncation method \cite{HLS2,HLS3}. Let us choose two $C^\ii$ functions $\chi$ and $\xi$ with values in $[0,1]$ such that $\chi^2+\xi^2=1\,$, $\chi$ has its support in $B(0,2)$ and $\chi\equiv1$ on $B(0,1)\,$. We denote $\chi_R(x):=\chi(x/R)$ and $\xi_R(x):=\xi(x/R)\,$. By the IMS localization formula, we have
$$-\Delta=\chi_R\,(-\Delta)\,\chi_R+\xi_R\,(-\Delta)\,\xi_R-|\nabla\chi_R|^2-|\nabla\xi_R|^2\;,$$
from which we deduce that
$$\tr\big((-\Delta)\,\gamma_n\big)\geq \tr\big((-\Delta)\,\chi_R\,\gamma_n\,\chi_R\big)+\tr\big((-\Delta)\,\xi_R\,\gamma_n\,\xi_R\big)- \frac{C\,q}{R^2}\;.$$
As $\rho_{\gamma_n}=\rho_{\chi_R\gamma_n\chi_R}+\rho_{\xi_R\gamma_n\xi_R}\,$, we have
\begin{equation*}
D(\rho_{\gamma_n},\rho_{\gamma_n}) \geq D(\rho_{\chi_R\gamma_n\chi_R},\rho_{\chi_R\gamma_n\chi_R})+D(\rho_{\xi_R\gamma_n\xi_R},\rho_{\xi_R\gamma_n\xi_R})
\end{equation*}
where we have used that
\begin{equation*}
D(\rho_{\chi_R\gamma_n\chi_R},\rho_{\xi_R\gamma_n\xi_R}) =\iint_{\R^6}\frac{\rho_{\gamma_n}(x)\,\rho_{\gamma_n}(y)\,\chi_R(x)^2\,\xi_R(y)^2}{|x-y|}\;dx\,dy\geq0\;.
\end{equation*}
Using the estimates $\chi_R(x)^2\,\xi_R(y)^2/|x-y|\leq 1/R$ if $|y|\geq3\,R\,$, and $\xi_R(y)^2\leq\xi_{3R}(y)^2$ if $|y|\leq3\,R\,$, we obtain
\begin{multline*}
\iint_{\R^6}\frac{|\gamma_n(x,y)|^2\,\chi_R(x)^2\,\xi_R(y)^2}{|x-y|}\;dx\,dy \leq \frac{1}{R}\iint_{\R^6}|\gamma_n(x,y)|^2\;dx\,dy\\
+ \iint_{\R^6}\frac{|\gamma_n(x,y)|^2\,\chi_R(x)^2\,\chi_{3R}(y)^2}{|x-y|}\;dx\,dy\;.
\end{multline*}
We may also observe that
$$\iint_{\R^6}|\gamma_n(x,y)|^2\;dx\,dy=\tr(\gamma^2)\leq \tr(\gamma)=q\;.$$
The last ingredient is the 
%---------------------------------------------------------------------------
\begin{lemma}[Brown-Kosaki's inequality \cite{BK} ]\label{lem:entropy_subadditive} Assume that $\beta$ satisfies \textbf{(A1)--(A3)}. Let $\gamma\in\cK$ and consider a self-adjoint operator $X:L^2(\R^3)\to L^2(\R^3)$ such that $X^2\leq 1\,$. Then
\begin{equation}
\tr\big(\beta(X\,\gamma\,X)\big)\leq \tr\big(X\,\beta(\gamma)\,X\big)\;.
\label{entropy_subadditive}
\end{equation} \end{lemma}
%---------------------------------------------------------------------------

Using \eqref{entropy_subadditive}, we obtain
$$\tr\big(\beta(\chi_R\,\gamma\,\chi_R)\big)+\tr\big(\beta(\xi_R\,\gamma\,\xi_R)\big)\leq \tr\big(\beta(\gamma)\big)\;.$$
Hence we obtain
$$\cE_Z^\beta(\gamma_n) \geq \cE_Z^\beta(\chi_R\,\gamma_n\,\chi_R)+\cE_0^\beta(\xi_R\,\gamma_n\,\xi_R)-\frac{C\,q}{R} - \iint_{\substack{R\leq |x|\leq 2R\\R\leq |y|\leq 3R}}\frac{|\gamma_n(x,y)|^2}{|x-y|}\;dx\,dy\;.$$
Notice that $\tr\big( \xi_R\,\gamma_n\,\xi_R\big)=q-\,\tr\big( \chi_R\,\gamma_n\,\chi_R\big)$ and $I^\beta_0\big(q-\,\tr\big( \chi_R\,\gamma_n\,\chi_R\big)\big)=0$ by Lemma~\ref {lem:value_free}. Summarizing, we have proved that
\begin{equation}
\cE_Z^\beta(\gamma_n) \geq \cE_Z^\beta(\chi_R\,\gamma_n\,\chi_R)-\frac{C\,q}{R} - \iint_{\substack{R\leq |x|\leq 2R\\R\leq |y|\leq 3R}}\frac{|\gamma_n(x,y)|^2}{|x-y|}\;dx\,dy\;.
\label{estim_energy_splitted}
\end{equation}%---------------------------------------------------------------------------
\begin{lemma}\label{lem:compacite_gamma} We have $\chi_R\,\gamma_n\,\chi_R\to \chi_R\,\gamma\,\chi_R$ strongly in $\gS_1\,$, as $n\to\ii\,$. Hence
$$\lim_{n\to\ii}\tr\big(\chi_R\,\gamma_n\,\chi_R\big)=\tr\big(\chi_R\,\gamma\,\chi_R\big)\;.$$
\end{lemma}
%---------------------------------------------------------------------------
\begin{proof} We have
$$\tr\big(\chi_R\,\gamma_n\,\chi_R\big) = \tr\big((1-\Delta)^{1/2}\,\gamma_n\,(1-\Delta)^{1/2}(1-\Delta)^{-1/2}\,\chi_R^2\,(1-\Delta)^{-1/2}\big)\;,$$
Since $(1-\Delta)^{1/2}\,\gamma_n\,(1-\Delta)^{1/2}\wto(1-\Delta)^{1/2}\,\gamma\,(1-\Delta)^{1/2}$ weakly in $\gH$ and $(1-\Delta)^{-1/2}\,\chi_R^2\,(1-\Delta)^{-1/2}$ is a compact operator for any fixed $R\,$, we get
\begin{align*}
\lim_{n\to\ii}\tr\big(\chi_R\,\gamma_n\,\chi_R\big)& =\tr\big((1-\Delta)^{1/2}\,\gamma\,(1-\Delta)^{1/2}(1-\Delta)^{-1/2}\,\chi_R^2\,(1-\Delta)^{-1/2}\big)\\
& = \tr\big(\chi_R\,\gamma\,\chi_R\big)\;.
\end{align*}\end{proof}

Passing to the limit in \eqref{estim_energy_splitted} as $n\to\ii$ using the value of $I^\beta_0$ and the weak lower semi-continuity of $\cE_Z^\beta$ as proved in Lemma~\ref{lem:weakly_continuous_and_coercive}, we obtain
\begin{equation*}
I^\beta_Z(q) \geq \cE_Z^\beta(\chi_R\,\gamma\,\chi_R) -\frac{C\,q}{R} - \iint_{\substack{R\leq |x|\leq 2R\\R\leq |y|\leq 3R}}\frac{|\gamma(x,y)|^2}{|x-y|}\;dx\,dy\;.
\label{estim_energy_splitted2}
\end{equation*}
Now we can pass to the limit as $R\to\ii$ and we obtain
\begin{equation*}
I^\beta_Z(q) \geq I^\beta_Z(q')\;,
\end{equation*}
which contradicts $(ii)\,$.

\paragraph*{Step 4. The self-consistent equation.}\hfill
%---------------------------------------------------------------------------
\begin{lemma} Assume that $\beta$ satisfies \textbf{(A1)--(A4)} with $T=1$ and consider a minimizer $\gamma$ for $I^\beta_Z(q)\,$. Then, for some $\mu\leq0\,$, $\gamma$ satisfies the self-consistent equation
$$\gamma=g\big(H_\gamma-\mu\big)\;.$$
\end{lemma}
%---------------------------------------------------------------------------
\begin{proof} The minimization problem is set on a convex set, so $\gamma$ solves the linearized problem
$$\inf_{\gamma'\in\cK_q} \left\{\tr\big((H_\gamma+\beta'(\gamma))\,\gamma'\big)\right\}\;.$$
The function $\beta$ being strictly convex, $\gamma$ also solves the convex minimization problem
\begin{equation*}
\inf_{\gamma'\in\cK_q} \left\{\tr\left(H_\gamma\,\gamma'+\beta(\gamma')\right)\right\}\;.
\label{linearized_pb}
\end{equation*}
This minimizer is unique and given by $g(H_\gamma-\mu)$ where the multiplier $\mu$ is chosen to ensure the constraint $\tr\gamma=q$ (this is an easy adaptation of \eqref{value_linear}). Hence
$$\gamma=g\big(H_\gamma-\mu\big)\;.$$
Notice that $\sigma_{\rm ess}\big(H_\gamma\big)=[0,\ii)\,$, hence necessarily $\mu\leq0\,$, since $\gamma$ is a trace-class operator and $g>0$ on $(-\infty,0)$. \end{proof}

\paragraph*{Step 5. Non-Existence if $q\geq 2\,Z+1\,$.} The proof of Lieb \cite{Lieb-84} for the usual Hartree-Fock case at zero temperature applies here. We only sketch it for the convenience of the reader.

Assume that $I^\beta_Z(q)$ admits a minimizer $\gamma\,$. It satisfies the self-consistent equation $\gamma=g(H_\gamma-\mu)\,$. In particular, $[\gamma,H_\gamma]=0$ and $H_\gamma\,\gamma\leq0\,$. Hence we have
$$\tr(|x|\,H_\gamma\,\gamma)\leq0\;.$$
Inserting the definition of $H_\gamma\,$, we obtain
\begin{equation}
0\geq \tr(|x|\,(-\Delta)\,\gamma)-Z\,q +\iint_{\R^6}\frac{|x|}{|x-y|}\,\big(\rho_\gamma(x)\,\rho_\gamma(y)-|\gamma(x,y)|^2\big)\;dx\,dy\;.
\label{estim_Lieb}
\end{equation}
Now we have by symmetry
\begin{align*}
&\iint_{\R^6}\frac{|x|}{|x-y|}\,\big(\rho_\gamma(x)\,\rho_\gamma(y)-|\gamma(x,y)|^2\big)\;dx\,dy\\
&\qquad=\frac12 \iint_{\R^6}\frac{(|x|+|y|)}{|x-y|}\,\big(\rho_\gamma(x)\,\rho_\gamma(y)-|\gamma(x,y)|^2\big)\;dx\,dy\\
&\qquad\geq\frac12\iint_{\R^6}\big(\rho_\gamma(x)\,\rho_\gamma(y)-|\gamma(x,y)|^2\big)\;dx\,dy=\frac12\big(q^2-\,\tr(\gamma^2)\big)\geq \frac12\big(q^2-q\big)
\end{align*}
where we have used the fact that $\tr\gamma^2\leq\tr\gamma=q$ because $0\leq\gamma\leq 1\,$. Using that $|x|\,(-\Delta)+(-\Delta)\,|x|\geq0$ which is equivalent to Hardy's inequality, we obtain $\tr\big(|x|\,(-\Delta)\,\gamma\big)>0\,$. Inserting in \eqref{estim_Lieb}, this yields $q< 2\,Z+1\,$.

\paragraph*{Step 6. Existence for $I^\beta_Z\,$.} By Lemma~\ref{lem:value_free} and Lemma~\ref{lem:binding_inequalities}, we have that $q\mapsto I^\beta_Z(q)$ is non-increasing. Condition $(ii)$ can be written $I^\beta_Z(q)<I^\beta_Z(q')$ for all $q'\in[0,q)\,$. As we know that for $q\geq 2\,Z+1\,$, $I^\beta_Z(q)$ has no minimizer, necessarily $I^\beta_Z$ is constant on $[2\,Z+1,\ii)$:
$$\forall\, q\geq2\,Z+1\;,\quad I^\beta_Z(q)=I^\beta_Z(2\,Z+1)\;.$$
Hence
$$I^\beta_Z=\inf\left\{\cE_Z^\beta(\gamma)\;|\;\gamma\in\cK\,,\;\tr\gamma\leq 2\,Z+1\right\}\;.$$
Now $\cE_Z^\beta$ is coercive and wlsc on $\{\gamma\in\cK\,|\;\tr\gamma\leq 2\,Z+1\}$ by Lemma~\ref{lem:weakly_continuous_and_coercive}, hence it admits a minimizer, with a vanishing Lagrange multiplier $\mu\,$.

\medskip This ends the proof of Theorem~\ref{thm_exists}.\hfill\qed

%%%%%%%%%%%%%%%%%%%%%%%%%%%%%%%%%%%%%%%%%%%%%%%%%%%%%%%%%%%%%%%%%%%%%%%%%%%%
\subsection{Proof of Proposition~\ref{prop_exists_g}}\label{sec:prop_exists_g}

We start the proof by two preliminary results.
%---------------------------------------------------------------------------
\begin{lemma}\label{I_negative} For any $q>0\,$, $I^\beta_Z(q)<0$ and $\lim_{q\to0_+}I^\beta_Z(q)=0=I^\beta_Z(0)\,$. \end{lemma}
%---------------------------------------------------------------------------
\begin{proof} By Lemma~\ref{lem:binding_inequalities}, $I^\beta_Z(q)$ is nonincreasing, so it is enough to prove the result for $q>0\,$, small. Let us construct a state $\gamma\in\cK_q$ with $q$ small such that $\cE^\beta_Z(\gamma)<0\,$. If~$\phi$ is a normalized eigenvector of $-\Delta-Z/|x|$ corresponding to the first eigenvalue $-Z^2/4\,$, we choose
$$\gamma=q\,|\phi\rangle\langle\phi|\;.$$
The result follows from
$$\cE_Z^\beta(\gamma)=-\frac{q\,Z^2}4+\beta(q)=-\frac{q\,Z^2}4+o(q)\quad\mbox{as }\;q\to0_+\;,$$
since the direct and exchange term cancel for a rank-one projector. The fact that $I^\beta_Z(0)=0$ is a consequence of $\cap_{q\geq0}\left\{\gamma\in\cK\,|\,\tr(\gamma)\leq q\right\}=\{0\}\,$. \end{proof}

%---------------------------------------------------------------------------
\begin{lemma}\label{lem_estimate_eigenvalue_radial} Let $\gamma\in\cK_{q'}$ and denote by $\{\lambda_j\}_{j\geq1}$ the ordered sequence of negative eigenvalues of~$H_\gamma\,$. Then we have
\begin{equation}
\lambda_j\leq -\frac{(Z-q')^2}{4\,j^2}\;.
\label{estimate_eigenvalue_radial}
\end{equation}
In particular $H_\gamma$ has infinitely many negative eigenvalues converging to zero. \end{lemma}
%---------------------------------------------------------------------------
\begin{proof} The proof is inspired by several former papers \cite{Lieb-Simon-77,Lions-87}. We consider the subspace $V$ of $L^2(\R^3)$ consisting of radial functions. For any nonnegative function $\rho\in V$, we can estimate
$$\int_{\R^3}\frac{\rho(y)}{|x-y|}\;dy\leq \frac{\int_{\R^3}\rho\,dx}{|x|}$$
by Newton's Theorem. Recall that this is an easy consequence of the formula
$$u(x):=\int_{\R^3}\frac{\rho(y)}{|x-y|}\;dy=\frac 1{|x|}\int_{|x|<|y|}\rho\;dy+\int_{|x|>|y|}\frac{\rho(y)}{|y|}\;dy\;.$$
By differentiating with respect to $r=|x|$, it is indeed easy to check that $u$ is the unique radial solution to $(r^2\,u')'=-\,4\,\pi\,r^2$ such that $\lim_{r\to\infty}u(r)=0$. Hence, denoting $d\lambda(R)$ the Haar measure of $\mathrm{SO}_3\,$, if $\phi$ is a radial test function, we obtain
\begin{eqnarray*}
\iint_{\R^6}\frac{\phi^2(x)\,\rho_\gamma(y)}{|x-y|}\;dx\,dy &=&\int_{\mathrm{SO}_3}d\lambda(R)\iint_{\R^6}\frac{\phi^2(x)\,\rho_\gamma(R\,y)}{|x-y|}\;dx\,dy\\
& \leq & \int_{\R^3}\frac{\phi^2(x)\int_{\R^3}dy\int_{\mathrm{SO}_3}d\lambda(R)\,\rho_\gamma(R\,y)}{|x|}\;dx \\
&&\quad = \int_{\R^3}\frac{\phi^2(x)\,q'}{|x|}\;dx
\end{eqnarray*}
where Newton's Theorem has been applied to $\int_{\mathrm{SO}_3}\rho_\gamma(R\,y)\,d\lambda(R)\,$. This proves that
$$\Pi_V\,H_\gamma\,\Pi_V \leq \Pi_V\left({\textstyle -\Delta-\frac{Z-q'}{|x|}}\right)\Pi_V$$
where $\Pi_V$ is the orthogonal projector onto $V$. The eigenvalues of the operator of the right hand side, restricted to $V$, are non degenerate and equal to $-(Z-q')^2/(4\,j^2)\,$. By the theory of Weyl \cite{Reed-Simon}, we know that the essential spectrum of $H_\gamma$ is $[0,\infty)$, and so eigenvalues can accumulate only at $0_-\,$. \end{proof}

Let us prove Proposition~\ref{prop_exists_g}. We argue by contradiction and consider a $q\in(0,Z]$ satisfying \eqref{estim_q_exists} and for which $(ii)$ in Theorem~\ref{thm_exists} is false. This means that there exists $q'\in[0,q)$ such that $I^\beta_Z(q)=I^\beta_Z(q')\,$. By Lemma~\ref{lem:binding_inequalities}, $I^\beta_Z$ is constant on $[q',q]\,$. We take $q'$ as the smallest number satisfying this property. We notice that $q'$ satisfies \eqref{estim_q_exists} because the r.h.s. of \eqref{estim_q_exists} is nondecreasing. By Lemma~\ref{I_negative}, we have $q'>0\,$. Also $I^\beta_Z(q')<I^\beta_Z(q'')$ for any $q''\in[0,q')\,$. Hence $q'$ satisfies $(ii)$ and $I^\beta_Z(q')$ admits a minimizer $\gamma\,$. It satisfies the self-consistent equation
\begin{equation*}
\gamma=g(H_\gamma-\mu)
\label{SCF_absurde}
\end{equation*}
for some $\mu\leq0\,$.

If $\mu<0\,$, by Lemma~\ref{lem_estimate_eigenvalue_radial}, we can choose an eigenvalue $\lambda_j\in(\mu,0)$ with associated eigenfunction $\phi\,$. Notice that $g(\lambda_j-\mu)=g'(\lambda_j-\mu)=0$ and
$$0\leq \gamma+\nu\,|\phi\rangle\langle\phi| \leq 1$$
for any $\nu\in[0,1]\,$. Now we compute
\begin{equation*}
\cE_Z^\beta\big(\gamma+\nu\,|\phi\rangle\langle\phi|\big) = \cE_Z^\beta(\gamma)+\nu\,\lambda_j+\beta(\nu) = I_Z^\beta(q')+\nu\,\lambda_j+o(\nu)\;,
\end{equation*}
by Assumptions \textbf{(A1)--(A3)}. Hence $I^\beta_Z(q'+\nu)<I^\beta_Z(q')$ for any $\nu\in(0,q-q')\,$, a contradiction.

If $\mu=0$, then
$$q'=\tr\big(g(H_\gamma)\big)\geq \sum_{j\geq1}\, g\left({\textstyle\frac{-(Z-q')^2}{4\,j^2}}\right)\,,$$
from which we deduce that $q'=q\,$, again a contradiction.\hfill\qed

%%%%%%%%%%%%%%%%%%%%%%%%%%%%%%%%%%%%%%%%%%%%%%%%%%%%%%%%%%%%%%%%%%%%%%%%%%%%
%%%%%%%%%%%%%%%%%%%%%%%%%%%%%%%%%%%%%%%%%%%%%%%%%%%%%%%%%%%%%%%%%%%%%%%%%%%%
\subsection*{Concluding remarks}

The \HFT~is a generalization of the usual \HF, which corresponds to the zero temperature case. Many variants of this model can be proposed in the case of the euclidean space, each of them giving rise to minimizers which are in general finite rank operators, the rank being usually larger than the number of electrons. Generically we obtain mixed states.

Almost for free, an orbital stability result follows from the minimization scheme. Proving the stability with respect to a determined stationary state is still an open question. This would be possible if we knew that the minimizers are isolated, a property which is not granted for the zero temperature case and even known to be false in some special cases.

Our framework is very natural when dealing with density operators, of finite or infinite rank. Recovering know results for the standard \HF~by letting the temperature go to zero is not really difficult. This should not suggest that the non-zero temperature case is similar to the case with zero temperature. For instance, the maximal number of electrons or the total charge $q^{\rm HF}_{\rm max}(T)$ that can be binded depends on $T\,$, and the study of this ionization threshold as a function of~$T$ is an open question.

%%%%%%%%%%%%%%%%%%%%%%%%%%%%%%%%%%%%%%%%%%%%%%%%%%%%%%%%%%%%%%%%%%%%%%%%%%%%
%%%%%%%%%%%%%%%%%%%%%%%%%%%%%%%%%%%%%%%%%%%%%%%%%%%%%%%%%%%%%%%%%%%%%%%%%%%%
\bigskip\noindent{\it Acknowledgements.} {\small Authors have been supported by the ANR Accquarel and the ECOS project no. C05E09. P.F.  was  partially supported by Fondecyt Grant \# 1070314  and FONDAP de Matem\'aticas Aplicadas.}

\medskip\noindent{\scriptsize\copyright~2008 by the authors. This paper may be reproduced, in its entirety, for non-commercial~purposes.}

%%%%%%%%%%%%%%%%%%%%%%%%%%%%%%%%%%%%%%%%%%%%%%%%%%%%%%%%%%%%%%%%%%%%%%%%%%%%
%%%%%%%%%%%%%%%%%%%%%%%%%%%%%%%%%%%%%%%%%%%%%%%%%%%%%%%%%%%%%%%%%%%%%%%%%%%%
% \nocite*
% \bibliographystyle{siam}\small
% \bibliography{References}
% 
% \end{document}
%%%%%%%%%%%%%%%%%%%%%%%%%%%%%%%%%%%%%%%%%%%%%%%%%%%%%%%%%%%%%%%%%%%%%%%%%%%%
%%%%%%%%%%%%%%%%%%%%%%%%%%%%%%%%%%%%%%%%%%%%%%%%%%%%%%%%%%%%%%%%%%%%%%%%%%%%

\end{document}